\documentclass[10pt,conference,compsocconf]{IEEEtran}
\usepackage{flushend}
\usepackage{ifpdf}
\usepackage{cite}
\usepackage[cmex10]{amsmath}
\usepackage{algorithmic}
\usepackage{array}
\usepackage[font=footnotesize]{subfig}
\usepackage{fixltx2e}
\usepackage{url}
\usepackage[american]{babel}
\usepackage[utf8]{inputenc}
\usepackage{amsfonts,amssymb}
\usepackage{stmaryrd}
\usepackage{paralist}
\usepackage{xspace}
\usepackage[ruled,vlined]{algorithm2e}
\usepackage{graphicx,color,psfrag}
\graphicspath{{fig/}}
\usepackage{pstool}
\usepackage{multirow}
\usepackage{ifthen}
\usepackage{todonotes}
\usepackage{fullpage}

\newcommand{\GEMM}{\ensuremath{\mathit{GEMM}}\xspace}

\newcommand{\TRSM}{\ensuremath{\mathit{TRSM}}\xspace}

\newcommand{\HQR}{\textsc{HQR}\xspace}
\newcommand{\MC}{\textsc{Fibonacci}\xspace}

\newcommand{\Greedy}{\textsc{Greedy}\xspace}

\newcommand{\killer}{\ensuremath{\mathit{eliminator}}\xspace}

\newcommand{\parsec}{\textsc{PaRSEC}\xspace}

\newcommand{\factor}{\ensuremath{\mathit{factor}}\xspace}
\newcommand{\apply}{\ensuremath{\mathit{apply}}\xspace}
\newcommand{\killm}{\ensuremath{\mathit{eliminate}}\xspace}
\newcommand{\update}{\ensuremath{\mathit{update}}\xspace}
\newcommand{\luqralgo}{\textit{LU-QR~Algorithm}\xspace}
\newcommand{\dague}{PaRSEC\xspace}

\newcommand{\flop}{floating-point operations\xspace}
\newcommand{\glops}{FLOP/sec\xspace}

\title{Designing LU-QR hybrid solvers for performance and stability}

\author{\IEEEauthorblockN{
Mathieu Faverge$^{1}$, Julien Herrmann$^{2}$, Julien Langou$^{3}$,\\
Bradley Lowery$^{3}$, Yves Robert$^{2,4}$ and Jack Dongarra$^{4}$}
\IEEEauthorblockA{$1.$ Laboratoire LaBRI, IPB ENSEIRB-MatMeca, Bordeaux, France}
\IEEEauthorblockA{$2.$ Laboratoire LIP, \'Ecole Normale Sup\'erieure de Lyon, France}
\IEEEauthorblockA{$3.$ University Colorado Denver, USA}
\IEEEauthorblockA{$4.$ University of Tennessee Knoxville, USA}
}

\begin{document}

\maketitle

\begin{abstract}
This paper introduces hybrid LU-QR algorithms for solving dense linear systems of the form $Ax=b$. Throughout a matrix factorization,
these algorithms dynamically alternate LU with local pivoting and QR elimination steps, based upon some
robustness criterion. LU elimination steps can be very efficiently parallelized, 
and are twice as cheap in terms of \flop, as QR steps.
However, LU steps are not necessarily stable, while QR steps are always stable.
The hybrid algorithms execute a QR step when a robustness criterion detects some risk for instability,
and they execute an LU step otherwise. Ideally, the choice between LU and QR steps must have a small
computational overhead and must provide a satisfactory level of
stability with as few QR steps as possible.
In this paper, we introduce several robustness criteria and we establish
upper bounds on the growth factor of the norm of the updated matrix incurred by each of these criteria.
In addition, we describe the implementation of the hybrid algorithms through an extension of the \dague software 
to allow for dynamic choices during execution. Finally, 
we analyze both stability and performance results compared to state-of-the-art linear solvers
on parallel distributed multicore platforms.
\end{abstract}

\IEEEpeerreviewmaketitle

\section{Introduction}
\label{sec.intro}

Consider a dense linear system $A x = b$ to solve, where $A$ 
is a square tiled-matrix, with  $n$ tiles per row or column.
Each tile is a block of  $n_b$-by-$n_b$ elements, 
so that the actual size of $A$ is $N = n \times n_{b}$.
Here, $n_b$ is a
parameter tuned to squeeze the most out of arithmetic units and memory
hierarchy.
To solve the linear system $Ax =b$, with $A$ a general matrix, one usually applies a series of transformations, pre-multiplying $A$ by several elementary matrices. There are two main approaches:
\emph{LU factorization}, where one uses lower unit triangular matrices, and \emph{QR factorization}, where one uses orthogonal Householder matrices. 
To the best of our knowledge, this paper is the first study to propose a mix of both approaches during a single factorization.
The LU factorization update is based upon matrix-matrix multiplications, a kernel that can be very efficiently parallelized, and whose library implementations typically achieve close to peak CPU performance.
Unfortunately,  the  efficiency of LU factorization is hindered 
by the need to perform partial pivoting at each step of the algorithm,
to ensure numerical stability. On the contrary, the QR factorization is
always stable, but requires twice as many \flop, and a more
complicated update step that is not as parallel as a matrix-matrix
product. Tiled QR algorithms~\cite{Buttari2008,tileplasma,QuintanaOrti:2009} greatly improve the parallelism of the
update step since they involve no pivoting but are based upon
more complicated kernels whose library implementations requires twice as many operations as LU.

The main objective of this paper is to explore the design of  \emph{hybrid}
algorithms that would combine the low cost and high CPU efficiency of the
LU factorization, while retaining the numerical stability of the QR
approach. 
In a nutshell, the idea is the following: at each step of the elimination,
we perform a 
\emph{robustness} test to know if the diagonal tile can be stably used to eliminate the tiles beneath it using an LU step.
If the test succeeds, then go for an elimination step based upon 
LU kernels, without any further pivoting involving sub-diagonal tiles 
in the panel. Technically, this is very similar to a step during a block LU factorization~\cite{blocklu:95}.
Otherwise, if the test fails, then go for a step with QR kernels.
On the one extreme, if all tests succeed throughout the algorithm, we implement an LU factorization without
pivoting.  
On the other extreme, if all tests fail, we implement a QR factorization. On the average, only a fraction of the tests will fail. If this fraction remains small enough,
we will reach a CPU performance close to that of LU without pivoting.
Of course the challenge is to design a test that is accurate enough (and not too costly)
so that LU kernels are applied only when it is numerically safe to do so.

Implementing such a hybrid algorithm on a state-of-the-art  distributed-memory platform, whose nodes are themselves equipped with multiple cores, is a 
programming challenge. Within a node, 
the architecture is a shared-memory machine, running many parallel
threads on the cores.  But the global architecture is
a distributed-memory machine, and requires MPI communication primitives for
inter-node communications.  A slight change in the algorithm, or in the
matrix layout across the nodes, might call for a time-consuming and
error-prone process of code adaptation. For each version, one must  identify,
and adequately implement, inter-node versus intra-node kernels. This
dramatically complicates the task of the programmers if they rely on a manual
approach. We solve this problem by relying on 
the \dague software~\cite{bosilca2012dague,dague-engine, dague-la}, 
so that we can concentrate on the algorithm and forget about MPI and threads. Once we have specified the algorithm at a task level, the
\dague software will recognize which operations are local to a node (and
hence correspond to shared-memory accesses), and which are not (and hence must
be converted into MPI communications).  Previous experiments show that this approach is very powerful, and that the use of a higher-level framework does not prevent our algorithms from achieving the same performance 
as state-of-the-art library releases~\cite{dongarra2013hierarchical}.

However, implementing a hybrid algorithm requires the programmer to implement a \emph{dynamic} task graph of the computation.
Indeed, the task graph of the hybrid factorization algorithm
is no longer known statically (contrarily to a standard LU or QR factorization). At each step of the elimination,
we use either LU-based or QR-based tasks, but not both. This requires the algorithm 
to dynamically fork upon the outcome of the robustness test, in order to apply the
selected kernels. The solution is to prepare a graph that includes both types of  tasks, namely LU and QR kernels,
to select the adequate tasks on the fly, and to discard the useless ones. We have to join both potential
execution flows at the end of each step, symmetrically. Most of this mechanism is transparent to the user. We discuss this extension of \dague in more detail in Section~\ref{sec.dague}.

The major contributions of this paper are the following:
\begin{compactitem}
\item The introduction of new LU-QR hybrid algorithms;
\item The design of several robustness criteria, with bounds on the induced growth factor;
\item A comprehensive experimental evaluation of the best trade-offs between performance and numerical stability;
\item The extension of \dague to deal with dynamic task graphs.
\end{compactitem}
The rest of the paper is organized as follows. First we explain the main principles of LU-QR hybrid algorithms in Section~\ref{sec.luqr-algo}. Then we describe
robustness criteria in Section~\ref{sec.criteria}. Next we detail the implementation within the \dague framework  in Section~\ref{sec.dague}. We report
experimental results in Section~\ref{sec.experiments}. We discuss 
related work in Section~\ref{sec.related}. 
Finally, we provide concluding remarks and future directions in
Section~\ref{sec.conclusion}.

\section{Hybrid LU-QR algorithms}
\label{sec.luqr-algo}

In this section, we describe hybrid algorithms to solve a dense linear system
$A x = b$, where $A~=~(A_{ij})_{(i,j)\in\llbracket1..n\rrbracket^2}$ is a square tiled-matrix, with $n$
tiles per row or column.  Each tile is a block of  $n_b$-by-$n_b$ elements, so
that $A$ is of order $N = n \times n_{b}$. 

The common goal of a classical one-sided factorization (LU or QR)
is to {\em triangularize} the matrix $A$ through a succession of elementary transformations.
Consider the first step of such an algorithm. We partition $A$ by block such that $A = \begin{pmatrix}
  A_{11}    &  C  \\
   B   &  D
\end{pmatrix}$. In terms of tile,
$A_{11}$ is $1$-by-$1$,
$B$ is $(n-1)$-by-$1$, 
$C$ is $1$-by-$(n-1)$, and
$D$ is $(n-1)$-by-$(n-1)$.
The first block-column $\begin{pmatrix}
  A_{11}  \\
   B   
\end{pmatrix}$ is the {\em panel} of the current step.

Traditional algorithms (LU or QR) perform the same type of transformation at each step. The key observation of this paper
is that any type of transformation (LU or QR) can be used for a given step independently of what was used for the previous steps.
The common framework of a step is the following:
\begin{equation}
{\scriptsize \begin{pmatrix}
  A_{11}    &  C  \\
   B   &  D
\end{pmatrix}  \Leftrightarrow
\begin{pmatrix}
  \factor   &  \apply  \\
   \killm   & \update
\end{pmatrix}   \Leftrightarrow
\begin{pmatrix}
  U_{11}   &  C' \\
   0   & D'
\end{pmatrix}.
}
\label{eq.fourmatrices}
\end{equation}
First,
$A_{11}$ is {\em factored} and transformed in the upper triangular matrix $U_{11}$.
Then, the transformation of the factorization of $A_{11}$ is {\em applied} to $C$.
Then $A_{11}$ is used to {\em eliminate} $B$. Finally $D$ is accordingly {\em updated}. 
Recursively factoring $D'$ with the same framework will complete the factorization to an upper triangular matrix. 

For each step, we have a choice for an LU step or a QR step.
The operation count for each kernel is given in Table \ref{tab.count}.

\begin{table}[htdp]
\begin{center}
\begin{tabular}{|l|ll|ll|}
\hline
& \multicolumn{2}{c|}{LU step, var A1} &  \multicolumn{2}{c|}{QR step}\\ \hline
\factor $A$ & $2/3$ & GETRF & $4/3$ & GEQRT \\  \hline
\killm  $B$  & $(n-1)$ & TRSM & $2(n-1)$ & TSQRT \\ \hline
\apply $C$ &  $(n-1)$ & TRSM & $2(n-1)$ & TSMQR \\  \hline 
\update $D$  & $2 (n-1)^{2}$  & GEMM & $4(n-1)^{2}$ & UNMQR \\ \hline
\end{tabular}
\end{center}
\caption{Computational cost of each kernel. The unit is $n_{b}^{3}$ \flop.}
\label{tab.count}
\end{table}%

Generally speaking, QR transformations
are twice as costly as their LU counterparts. The bulk of the computations take place in the update of the trailing matrix $D$.
This obviously favors LU  \update kernels. In addition, the LU {\em update} kernels are fully parallel and 
can be applied independently on the $(n-1)^{2}$ trailing tiles. 
Unfortunately, LU updates (using GEMM) are stable only when $\| A_{11}^{-1}\|^{-1}$ is larger than $\| B \|$
(see Section~\ref{sec.criteria}). If this is not the case,
we have to resort to QR kernels. Not only these are twice as costly, but they also suffer from enforcing more dependencies: all columns can 
still be processed (\apply and \update kernels) independently, but inside a column, the kernels must be applied in sequence.

The hybrid \luqralgo  uses the standard  2D block-cyclic distribution of
tiles along a virtual $p$-by-$q$ cluster grid. The 2D block-cyclic distribution
nicely balances the load across resources for both LU and QR steps. Thus at
step $k$ of the factorization, the panel is split into $p$ {\em domains} of
approximately $\frac{n-k+1}{p}$ tile rows. .
Domains will be associated with physical memory regions, typically a domain per node in a distributed memory platform.
Thus an important design goal is to minimize the number of communications across domains, because these correspond to 
nonlocal communications between nodes.
At each step $k$ of the factorization, the domain of the node owning the diagonal tile $A_{kk}$ is
called the \emph{diagonal domain}.

The hybrid \luqralgo applies LU kernels when it is numerically safe to do so,
and QR kernels otherwise.
Coming back to the first elimination step, the sequence of operations is described in Algorithm~\ref{alg.LUQR}.


\begin{algorithm}[!t]
{\normalsize
	\For {$k = 1$ to $n$}{
		{\em Factor:} Compute a factorization of the diagonal tile: either with LU partial pivoting or QR\;
		{\em Check:} Compute some robustness criteria (see Section \ref{sec.criteria}) involving only tiles $A_{i,k}$, where $k \leq i \leq n$, in the elimination panel\;
		{\em Apply, Eliminate, Update:}\\
		\eIf{the criterion succeeds}{
			 Perform an LU step\;
		}
		{
			Perform a QR step\;
		}
	}
}
\caption{Hybrid LU-QR algorithm}
\label{alg.LUQR}
\end{algorithm}


\subsection{LU step}
\label{sec:A1} 
\label{sec:diagonaldomain}

We assume that the criterion validates an LU step (see Section~\ref{sec.criteria}).
We describe the variant (A1) of an LU step given in 
Algorithm~\ref{alg.LU}.

\begin{algorithm}[!t]
{\normalsize
	{\em Factor:} $A_{kk} \leftarrow GETRF(A_{kk})$ \;
	\For {$i = k+1$ to $n$}{
		{\em Eliminate:} $A_{ik} \leftarrow TRSM(A_{kk},A_{ik})$\;
	} 
	\For {$j = k+1$ to $n$}{
		{\em Apply:} $A_{kj} \leftarrow SWPTRSM(A_{kk},A_{kj})$\;
	} 
	\For {$i = k+1$ to $n$}{
		\For {$j = k+1$ to $n$}{
			{\em Update:} $A_{ij} \leftarrow GEMM(A_{ik},A_{kj},A_{ij})$\;
		} 
	} 
}
\caption{Step $k$ of an LU step - var (A1)}
\label{alg.LU}
\end{algorithm}

The kernels for the LU step are the following: 
\begin{compactitem}

\item {\em Factor:} $A_{kk} \leftarrow GETRF(A_{kk})$ is an LU factorization
with partial pivoting: $P_{kk} A_{kk} = L_{kk} U_{kk}$, the output matrices
$L_{kk}$ and $U_{kk}$ are stored in place of the input $A_{kk}$.

\item {\em Eliminate:} $A_{ik} \leftarrow TRSM(A_{kk},A_{ik})$ solves in-place, the 
upper triangular system such that $A_{ik} \leftarrow A_{ik}U_{kk}^{-1}$
where $U_{kk}$ is stored in the upper part of $A_{kk}$.

\item {\em Apply:} $A_{kj} \leftarrow SWPTRSM(A_{kk},A_{ik})$ 
solves the unit lower triangular system
such that $A_{kj} \leftarrow L_{kk}^{-1} P_{kk} A_{kj}$
where $L_{kk}$ is stored in the (strictly) lower part of $A_{kk}$.

\item {\em Update:} $A_{ij} \leftarrow GEMM(A_{ik},A_{kj},A_{ij})$ is a general matrix-matrix
multiplication $A_{ij} \leftarrow A_{ij}- A_{ik} A_{kj}$. 

\end{compactitem}

In terms of parallelism, 
the factorization of the diagonal tile is followed by the \TRSM kernels
that can be processed in parallel, then every \GEMM kernel can be processed concurrently. These 
highly parallelizable updates constitute one of the two main advantages of the LU step over the QR step.
The second main advantage is halving the number of \flop.

During the {\em factor} step, one variant is to factor the whole diagonal
domain instead of only factoring the diagonal tile. Considering
Algorithm~\ref{alg.LU}, the difference lies in the first line: rather than
calling $GETRF(A_{kk})$, thereby searching for pivots only within the diagonal
tile $A_{kk}$, we implemented a variant where we extend the search for pivots
across the \emph{diagonal domain} (the {\em Apply} step is modified
accordingly). Working on the diagonal domain instead of the diagonal tile
increases the smallest singular value of the factored region and therefore
increases the likelihood of an LU step. Since all tiles in the diagonal domain are local to a single node, 
extending the search to the diagonal domain is
done without any inter-domain communication.  The stability analysis of
Section~\ref{sec.criteria} applies to both scenarios, the one where $A_{kk}$ is
factored in isolation, and the one where it is factored with the help of the
diagonal domain.  In the experimental section, we will use the variant which
factors the diagonal domain.

\subsection{QR step}

If the decision to process a QR step is taken by the criterion, the LU decomposition of the diagonal 
domain is dropped, and the factorization of the panel starts over. This step of the factorization is 
then processed using orthogonal transformations. Every tile below the diagonal (matrix $B$ in Equation~\eqref{eq.fourmatrices})
is zeroed out using a 
triangular tile, or eliminator tile. In a QR step, the diagonal tile is 
factored (with a GEQRF kernel) and used to eliminate all the other tiles of the panel (with a TSQRT 
kernel) The trailing submatrix is updated, respectively, with UNMQR and
TSMQR kernels.
To further increase the degree of parallelism of the algorithm, it is 
possible to use several eliminator tiles inside a panel, typically one (or more) per domain.
The only condition is that concurrent 
elimination operations must involve disjoint tile pairs (the unique eliminator of tile $A_{ik}$ will 
be referred to as $A_{\killer(i,k),k}$). Of course, in the end, there must remain only one non-zero 
tile on the panel diagonal, so that all eliminators except the diagonal tile must be eliminated 
later on (with a TTQRT kernel on the panel and TTMQR updates on the
trailing submatrix), using a reduction tree of arbitrary shape. This reduction tree will involve inter-domain communications.
In our hybrid LU-QR algorithm, the QR step is processed following an instance of the generic 
hierarchical QR factorization \HQR~\cite{dongarra2013hierarchical} described in Algorithms~\ref{alg.QR} 
and~\ref{alg.elim}.   

\begin{algorithm}[!t]
{\normalsize
	\For {$i = k+1$ to $n$}{
		$elim(i,\killer(i,k),k)$\;
	} 
}
\caption{Step $k$ of the \HQR factorization}
\label{alg.QR}
\end{algorithm}

\begin{algorithm}[!t]
{\normalsize
	(a) With TS kernels\\
	$A_{\killer(i,k),k} \leftarrow GEQRT(A_{\killer(i,k),k})$\;
	$A_{i,k}, A_{\killer(i,k),k} \leftarrow TSQRT(A_{i,k}, A_{\killer(i,k),k})$\;
	\For {$j = k+1$ to $n-1$}{
		$A_{\killer(i,k),j} \leftarrow UNMQR(A_{\killer(i,k),j}, A_{\killer(i,k),k}$\;
		$A_{i,j}, A_{\killer(i,k),j} \leftarrow TSMQR(A_{i,j}, A_{\killer(i,k),j}, A_{i,k})$\;
	}

	\vspace{6mm}	
	
	(b) With TT kernels\\
	$A_{\killer(i,k),k} \leftarrow GEQRT(A_{\killer(i,k),k})$\;
	$A_{i,k} \leftarrow GEQRT(A_{i,k})$\;
	\For {$j = k+1$ to $n-1$}{
		$A_{\killer(i,k),j} \leftarrow UNMQR(A_{\killer(i,k),j}, A_{\killer(i,k),k}$\;
		$A_{i,j} \leftarrow UNMQR(A_{i,j}, A_{i,k}$\;
	}	
	$A_{i,k}, A_{\killer(i,k),k} \leftarrow TTQRT(A_{i,k}, A_{\killer(i,k),k})$\;
	\For {$j = k+1$ to $n-1$}{
		$A_{i,j}, A_{\killer(i,k),j} \leftarrow TTMQR(A_{i,j}, A_{\killer(i,k),j}, A_{i,k})$\;
	}
}	
\caption{\small Elimination $elim(i,\killer(i,k),k)$}
\label{alg.elim}
\end{algorithm}

\vspace{3mm}

Each elimination $elim(i, \killer (i, k), k)$ consists of two sub-steps: first 
in column $k$, tile $(i, k)$ is zeroed out (or killed) by tile $(\killer (i, k), k)$; and 
in each following column $j > k$, tiles $(i, j)$ and $(\killer (i, k), j)$ are updated; 
all these updates are independent and can be triggered as soon as the elimination is completed.
The algorithm is entirely characterized by its elimination list, which is the ordered 
list of all the eliminations $elim(i, \killer (i, k), k)$ that are executed. The orthogonal
transformation $elim(i, \killer(i, k), k)$ uses either a TTQRT kernel or a TSQRT kernel 
depending upon whether the tile to eliminate is either triangular or square.
In our hybrid \luqralgo, any combination of reduction trees of the
\HQR algorithm described in~\cite{dongarra2013hierarchical} is available.
It is then possible to use an
intra-domain reduction tree to locally eliminate many tiles without
inter-domain communication. A unique triangular tile is left on each 
node and then the reductions across domains are performed
following a second level of reduction tree.

\subsection{LU step variants}

In the following, we describe several other variants of the LU step.

\subsubsection{Variant (A2)} It consists of first performing a $QR$ factorization of the diagonal
tile and proceeds pretty much as in (A1) thereafter.
\begin{compactitem}

\item {\em Factor:} $A_{kk} \leftarrow GEQRF(A_{kk})$ is a QR factorization
$ A_{kk} = Q_{kk} U_{kk}$,
where $Q_{kk}$ is never constructed explicitly and we instead store
the Householder reflector $V_{kk}$.
The output matrices
$V_{kk}$ and $U_{kk}$ are stored in place of the input $A_{kk}$.

\item {\em Eliminate:} $A_{ik} \leftarrow TRSM(A_{kk},A_{ik})$ solves in-place the 
upper triangular system such that $A_{ik} \leftarrow A_{ik}U_{kk}^{-1}$
where $U_{kk}$ is stored in the upper part of $A_{kk}$.

\item {\em Apply:} $A_{kj} \leftarrow ORMQR(A_{kk},A_{ik})$ 
performs $A_{kj} \leftarrow Q_{kk}^{T} A_{kj}$
where $Q_{kk}^T$ is applied using $V_{kk}$ stored in the (strictly) lower part of $A_{kk}$.

\item {\em Update:} $A_{ij} \leftarrow GEMM(A_{ik},A_{kj},A_{ij})$ is a general matrix-matrix
multiplication $A_{ij} \leftarrow A_{ij}- A_{ik} A_{kj}$. 

\end{compactitem}
The {\em Eliminate} and {\em Update} steps are the exact same as in (A1).
The (A2) variant has the same data dependencies as (A1) and therefore the same level of 
parallelism.  A benefit of (A2) over (A1) is that if the criterion test decides that the step is a QR step,  then the
factorization of $A_{kk}$ is not discarded but rather used to continue the QR step.
A drawback of (A2) is that the {\em Factor} and {\em Apply} steps are twice as expensive as the ones in (A1).

\subsubsection{Variants (B1) and (B2)} Another option is to use the so-called {\em block LU factorization}~\cite{blocklu:95}.
The result of this formulation is a factorization
where the $U$ factor is block upper triangular  (as opposed to upper triangular), and
the diagonal tiles of the
$L$ factor are identity tiles. 
The {\em Factor} step can either be done using an LU factorization (variant (B1)) or a QR factorization (variant (B2)).
The {\em Eliminate} step is $A_{ik} \leftarrow A_{ik}A_{kk}^{-1}$. There is no
{\em Apply} step. And the {\em Update} step is $A_{ij} \leftarrow A_{ij}- A_{ik} A_{kj}$. 

The fact that row $k$ is not updated provides two benefits: 
(i) $A_{kk}$ does not need to be broadcast to these tiles,  simplifying the communication pattern;
(ii)  The stability of the LU step can be determined by considering only the growth factor in
the Schur complement of $A_{kk}$.  
One drawback of (B1) and (B2) is that the final matrix is not upper triangular but only block upper triangular.
This complicates the use of these methods to solve a linear system of equations. The stability of (B1) and (B2) 
has been analyzed in~\cite{blocklu:95}.

\medskip
We note that (A2) and (B2) use a QR factorization during the {\em Factor} step. Yet, we still call this an LU step.
This is because
all four LU variants mentioned use the Schur complement to update the trailing sub-matrix. 
The mathematical operation is:
$A_{ij} \leftarrow A_{ij} - A_{ik} A_{kk}^{-1}A_{kj} ,$. In practice, the {\em Update} step for all four variants
looks like $ A_{ij} \leftarrow A_{ij} - A_{ik} A_{kj} $, since $A_{kk}^{-1}$ is somehow applied to $A_{ik} $ and $ A_{kj} $ during the preliminary
{\em update} and {\em eliminate} steps.
The Schur update dominates
the cost of an LU factorization and therefore all variants are more efficient than a QR step.  Also, we have the same level
of parallelism for the update step: embarrassingly parallel.
In terms of stability, all variants would follow closely the analysis of Section~\ref{sec.stability}.
We do not consider further variants (A2), (B1), and (B2) in this paper, since they are all very similar,
and only study Algorithm~\ref{alg.LU}, (A1).

\subsection{Comments}

\subsubsection{Solving systems of linear equations}

To solve systems of linear equations, 
we augment $A$ with the right-hand side $b$ to get $\tilde{A} = (A,b)$
and apply all transformations to $\tilde{A}$. Then an $N$-by-$N$ triangular solve is needed.
This is the approach we used in our experiments.
We note that,
at the end of the factorization,
all needed information about the transformations is 
stored in place of $A$, so, alternatively, one can apply the transformations on $b$ during a second pass.

\subsubsection{No restriction on $N$}

In practice, $N$ does not have to be a
multiple of $n_b$. We keep this restriction for the sake of simplicity.
The algorithm can accommodate any $N$ and $n_b$ with some clean-up codes, which we have written.

\section{Robustness criteria}
\label{sec.criteria}

The decision to process an LU or a QR step is done dynamically during the
factorization, and constitutes the heart of the algorithm. Indeed, the decision
criteria has to be able to detect a potentially ``large'' stability
deterioration (according to a threshold) due to an LU step before its actual
computation,  in order to preventively switch to a QR step. As explained in
Section \ref{sec.luqr-algo}, in our hybrid LU-QR algorithm, the diagonal
tile is factored using an LU decomposition with partial pivoting.
At the same time, some data (like the norm of non local tiles
belonging to other domains) are collected and exchanged (using a Bruck's
all-reduce algorithm~\cite{Bruck1997}) between all nodes hosting at least
one tile of the panel. Based upon this information, all nodes make the decision
to continue the LU factorization step or to drop the LU decomposition of the
diagonal tile and process a full QR factorization step. The decision is
broadcast to the other nodes not involved in the panel factorization
within the next data communication. The decision process cost will depend on
the choice of the criterion and must not imply a large computational overhead
compared to the factorization cost. A good criterion will detect only the
``worst'' steps and will provide a good stability result with as few QR steps
as possible. In this section, we present three criteria, going from the most
elaborate (but also most costly) to the simplest ones. 

The stability of a step is determined by the growth of the norm of the updated matrix.
If a criterion determines the potential for an unacceptable growth due to an LU step, then
a QR step is used. A QR step is stable as there is no growth in the norm (2-norm) 
since it is a unitary transformation.
Each criterion depends on a 
threshold $\alpha$ that allows us to tighten or loosen the stability
requirement, and thus influence 
the amount of LU steps that we can afford during the factorization.  The optimal choice of $\alpha$
is not known. In Section~\ref{sec.stability}, we experiment with different choices of $\alpha$ for each
criterion.

\subsection{Max criterion}
LU factorization with partial pivoting chooses the largest element of a column as the pivot element.  
Partial pivoting is accepted as being numerically stable. 
However, pivoting across nodes is expensive.  To avoid this pivoting, we generalize the criteria to tiles
and determine if the diagonal tile is an acceptable pivot.   
A step is an LU step if
\begin{equation}
	\alpha \times \|( A_{kk}^{(k)} )^{-1} \|_1^{-1} \ge \max_{i > k} \|A_{i,k}^{(k)}\|_1.
	\label{eq:ppc}
\end{equation}
For the analysis we do not make an assumption as
to how the diagonal tile is factored.  We only assume that the diagonal tile is factored in a stable way 
(LU with partial pivoting or QR are acceptable). 
Note that, for the variant using pivoting in the diagonal domain
(see Section~\ref{sec:diagonaldomain}), which is the variant we experiment with in 
Section~\ref{sec.experiments},
$A_{kk}^{(k)}$ represents the
diagonal tile after pivoting among tiles in the diagonal domain.

To assess the growth of the norm of the updated matrix,
consider the update of the trailing sub-matrix.
For all $i,j > k$ we have:
\begin{align*}
	\|A_{i,j}^{(k+1)}\|_1 
		&=		\|A_{i,j}^{(k)} - A_{i,k}^{(k)} ( A_{k,k}^{(k)} )^{-1} A_{k,j}^{(k)} \|_1 \\
		&\le	\|A_{i,j}^{(k)}\|_1 +\| A_{i,k}^{(k)}\|_1\| ( A_{k,k}^{(k)} )^{-1} \|_1   \| A_{k,j}^{(k)} \|_1  \\
		&\le	\|A_{i,j}^{(k)}\|_1 + \alpha \| A_{k,j}^{(k)} \|_1 \\
		&\le	(1 + \alpha) \max\left( \|A_{i,j}^{(k)}\|_1, \| A_{k,j}^{(k)} \|_1 \right) \\
        	&\le	(1 + \alpha) \max_{i\ge k} \left( \|A_{i,j}^{(k)}\|_1 \right).
\end{align*}
The growth of any tile in the trailing sub-matrix is bounded by $1+\alpha$ times the largest tile in the same column.
If every step satisfies~\eqref{eq:ppc}, then we have the following bound:
\[ \frac{ \max_{i,j,k} \| A_{ij}^{(k)} \|_1 }{ \max_{i,j} \| A_{i,j} \|_1 } \le (1+\alpha)^{n-1}. \]
The expression above is a growth factor on the norm of the tiles.  For $\alpha = 1$, the growth factor of $2^{n-1}$ is an analogous result to an LU factorization with 
partial pivoting (scalar case)~\cite{Higham2002}. 
Finally, note that we can obtain this bound by generalizing the standard example for partial pivoting.  
The following matrix will match the bound above:
\[ A = \left( \begin{array}{cccc}
	\alpha^{-1} & 0 & 0 & 1 \\
	-1 & \alpha^{-1} & 0 & 1 \\
	-1 & -1 & \alpha^{-1} & 1 \\
	-1 & -1 & -1 & 1 
	\end{array}
	\right).
\]

\subsection{Sum criterion}
A stricter criteria is to compare the diagonal tile to the sum of the off-diagonal tiles:
\begin{equation}\label{eq:ddc}
\alpha \times \|( A_{kk}^{(k)} )^{-1} \|_{{1}}^{-1} \ge \sum_{i > k} \|A_{i,k}^{(k)}\|_1.
\end{equation}
Again, for the analysis, we only assume $A_{kk}^{-1}$ factored in a stable way.
For $\alpha \ge 1$, this criterion (and the Max criterion) is satisfied at every step if $A$ is block diagonally dominant~\cite{Higham2002}.  That is, a general matrix $A \in \mathbb{R}^{n \times n}$ 
is block diagonally dominant by columns with respect to a given partitioning $A = (A_{ij})$ and a given norm $\|\cdot\|$ if:
\[ \forall j \in \llbracket 1,n \rrbracket, \|A_{jj}^{-1}\|^{-1} \geq \sum_{i \neq j} \|A_{ij}\|. \]
Again we need to evaluate the growth of the norm of the updated trailing sub-matrix. For all $i,j > k$,
we have 
\begin{align*}
	\sum_{i > k} \|A_{i,j}^{(k+1)}\|_1 
		&=		\sum_{i > k} \|A_{i,j}^{(k)} - A_{i,k}^{(k)} ( A_{k,k}^{(k)} )^{-1} A_{k,j}^{(k)} \|_1 \\
		&\le	\sum_{i > k} \|A_{i,j}^{(k)}\|_1 \\
		&\quad+ \| A_{k,j}^{(k)} \|_1 \| ( A_{k,k}^{(k)} )^{-1} \|_1 \sum_{i > k} \| A_{i,k}^{(k)}\|_1 \\
		&\le	\sum_{i > k} \|A_{i,j}^{(k)}\|_1 + \alpha \| A_{k,j}^{(k)} \|_1.\\
\end{align*}
Hence, the growth of the updated matrix can be bounded in terms of an entire column rather than just an
individual tile.  The only growth in the sum is due to the norm of a single tile.   
For $\alpha =1$, the inequality becomes 
\[ \sum_{i > k} \|A_{i,j}^{(k+1)}\|_1 \le	\sum_{i \ge k} \|A_{i,j}^{(k)}\|_1. \]
If every step of the algorithm satisfies~\eqref{eq:ddc} (with $\alpha = 1$), then by induction we have:
\[ \sum_{i > k} \|A_{i,j}^{(k+1)}\|_1 \le \sum_{i \ge 1} \|A_{i,j}\|_1, \]
for all $i,j,k$.  This leads to the following bound:
\[ \frac{ \max_{i,j,k} \| A_{ij}^{(k)} \|_1 }{ \max_{i,j} \| A_{i,j} \|_1 } \le n. \]
From this we see that the criteria eliminates the potential for exponential growth due to the LU steps.
Note that for a diagonally dominant matrix, the bound on the growth factor can be reduced to 2~\cite{Higham2002}.

\subsection{MUMPS criterion}
In LU decomposition with partial pivoting, the largest element of the column is use as the pivot. This method is stable experimentally, but the seeking of the maximum and the pivoting requires a lot of communications in distributed memory. Thus in an LU step of the \luqralgo, the LU decomposition with partial pivoting is limited to the local tiles of the panel (i.e., to the diagonal domain). 
The idea behind the MUMPS criterion is to estimate the quality of the pivot found locally compared to the rest of the column. 
The MUMPS criterion is one of the strategies available in MUMPS although it is
for symmetric indefinite matrices ($LDL^T$)~\cite{Duff2005}, and  Amestoy et al.~\cite{mumps-crit} provided
us with their scalar criterion for the LU case.

At step $k$ of the \luqralgo, let $L^{(k)} U^{(k)}$ be the LU decomposition of the diagonal domain and $A_{ij}^{(k)}$ be the value of the tile $A_{ij}$ at the beginning of step $k$. Let $local\_max_k(j)$ be the largest element of the column $j$ of the panel in the diagonal domain, $away\_max_k(j)$ be the largest element of the column $j$ of the panel off the diagonal domain, and $pivot_k$ be the list of pivots used in the LU decomposition of the diagonal domain:
\begin{align*}
&local\_max_k(j) = \max_{\substack{\text{tiles } A_{i,k} \text{ on the} \\ \text{diagonal domain}}} ~\max_{l}~|(A_{i,k})_{l,j}|, \\
&away\_max_k(j) = \max_{\substack{\text{tiles } A_{i,k} \text{ off the} \\ \text{diagonal domain}}} ~\max_{l}~|(A_{i,k})_{l,j}|, \\
&pivot_k(j) = | U^{(k)}_{j,j} |.
\end{align*}
$pivot_k(j)$ represents the largest local element of the column $j$ at step $j$ of the LU decomposition with partial pivoting on the diagonal domain. Thus, we can express the growth factor of the largest local element of the column $j$ at step $j$ as: $growth\_factor_k(j) = pivot_k(j) / local\_max_k(j)$. The idea behind the MUMPS criterion is to estimate if the largest element outside the local domain would have grown the same way. Thus, we can define a vector $estimate\_max_k$ initialized to $away\_max_k$ and updated for each step $i$ of the LU decomposition with partial pivoting like $estimate\_max_k(j) \leftarrow estimate\_max_k(j) \times growth\_factor_k(i)$. We consider that the LU decomposition with partial pivoting of the diagonal domain can be used to eliminate the rest of the panel if and only if all pivots are larger than the estimated maximum of the column outside the diagonal domain times a threshold $\alpha$. Thus, the MUMPS criterion (as we implemented it) decides that step $k$ of the \luqralgo will be an LU step if and only if:
\begin{equation}\label{eq.mumps}
\forall j, \alpha \times pivot_k(j) \ge estimate\_max_k(j).
\end{equation}

\subsection{Complexity}
All criteria require the reduction of information of the off-diagonal tiles to the diagonal tile. Criteria~\eqref{eq:ppc} and~\eqref{eq:ddc} require the norm of
each tile to be calculated locally (our implementation uses the 1-norm) and then reduced to the diagonal tile.
Both criteria also require computing $\| A_{kk}^{-1} \|$.  Since the LU factorization of the 
diagonal tile is computed, the norm can be approximated using the L and U factors by an iterative method in $O(n_b^2)$ \flop. The overall complexity for both criteria is $O(n \times n_b^2 )$. 
Criterion~\eqref{eq.mumps} requires the maximum of each column be calculated locally and then reduced to the diagonal tile.  The complexity of the MUMPS criterion
is also $O(n \times n_b^2)$ comparisons. 

The Sum criterion  is the strictest of the three criteria.  It also provides the best  
stability with linear growth in the norm of the tiles in the worst case. The other two criteria have similar
worst case bounds.  The growth factor for both criteria are bound by the growth factor of partial (threshold) pivoting.
The Max criterion has a bound for the growth factor on the norm of the tiles that is 
analogous to partial pivoting.  The MUMPS criteria does not operate at the tile level, but rather on scalars.
If the estimated growth factor computed by the criteria is a good estimate, then the growth factor is
no worse than partial (threshold) pivoting. 

\section{Implementation}
\label{sec.dague}

As discussed in section~\ref{sec.intro}, we have implemented the \luqralgo on top of
the \parsec runtime. There are two major reasons for this choice: (i) it allows for
easily targeting distributed architectures while concentrating only on
the algorithm and not on implementation details such as data distribution
and communications; (ii) previous implementations of the \HQR
algorithm~\cite{dongarra2013hierarchical} can be reused for QR elimination steps, and they include efficient reduction trees
to reduce the critical path of these steps.

\begin{figure*}[ht]
\centering
\includegraphics[width=0.8\linewidth]{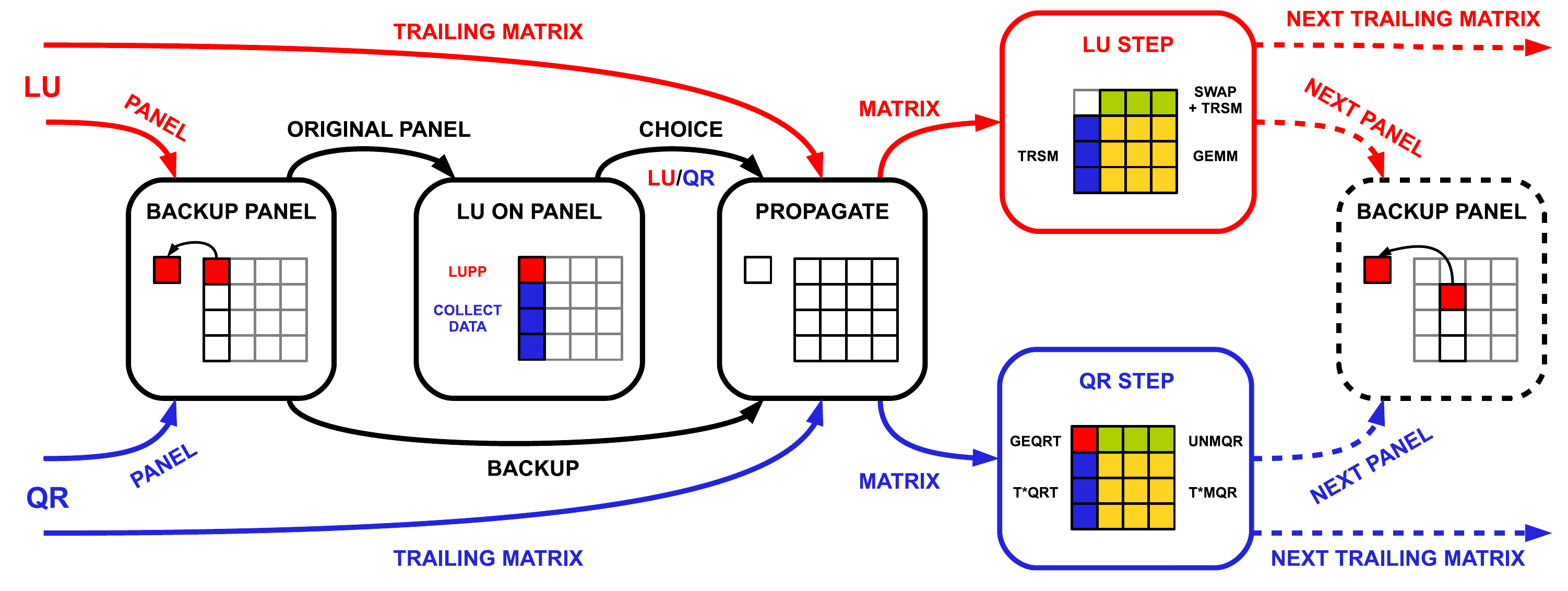}
\caption{\label{fig.choix}Dataflow of one step of the algorithm.}
\end{figure*}

However, this choice implied major difficulties due to the
parameterized task graph representation exploited by the \parsec
runtime. This representation being static, a solution had to be
developed to allow for dynamism in the graph traversal. To solve this
issue, a layer of selection tasks has been inserted between each
elimination step of the algorithm. These tasks are only executed
once a control flow has been sent to them after the criterion
selection. Thus, they delay the decision to send the data to the next
elimination step until a choice has been made, in order to guarantee that
data follow the correct path. These are the \emph{Propagate} tasks
on Figure~\ref{fig.choix}. Note that these tasks, as
well as  \emph{Backup Panel} tasks, can receive the same data from
two different paths. In the \parsec runtime, tasks are created only when
one of their dependencies is solved; then by graph construction they
are enabled only when the previous elimination step has already started, hence
they will receive their data only from the correct path. 

Figure~\ref{fig.choix} describes the connection between the
different stages of one elimination step of the algorithm. These
stages are described below:
\paragraph*{\textsc{Backup Panel}} This is a set of tasks that
  collect the tiles of the panel from the previous step. Since an LU
  factorization will be performed in-place for criterion computation,
  a backup of the tiles on the diagonal domain is created and 
  directly sent to the \emph{Propagate} tasks in case a QR elimination step
  is needed. All tiles belonging to the panel are forwarded to the
  \emph{LU On Panel} tasks.

\paragraph*{\textsc{LU On Panel}} Once the backup is done, the criterion
  is computed. The first node computes the $U$ matrix
  related to this elimination step. This could be done through LU
  factorization with or without pivoting. We decided to exploit the
  multi-threaded recursive-LU kernel from the PLASMA library to
  enlarge the pivot search space while keeping good efficiency~\cite{CPE:CPE3110}.
  All other nodes compute information required for the criterion (see
  section~\ref{sec.criteria}). Then, an all-reduce operation is
  performed to exchange the information, so that
  everyone can take and store the decision in a local array.  Once the
  decision is known, data on panel are forwarded to the appropriate
  \emph{Propagate} tasks and a control flow triggers all to release
  the correct path in the dataflow.

\paragraph*{\textsc{Propagate}} These tasks, one per tile, receive the
  decision from the previous stage through a control flow, and are
  responsible for forwarding the data to the computational tasks of
  the selected factorization. 
  The tasks belonging to the panel (assigned to the first nodes) have to restore the data back to
  their previous state if QR elimination is chosen. In all cases, the
  backup is destroyed upon exit of these tasks.

We are now ready to complete the description of each step:

\paragraph*{a) \textsc{LU Step}} If the numerical criterion is met by the panel
  computation, the update step is performed. On the nodes with the
  \emph{diagonal} row, the update is made according to the
  factorization used on the panel. Here, a swap is performed with all
  tiles of the local panel, and then a triangular solve is applied to the
  first row. On all other nodes, a block LU algorithm is used to
  performed the update. This means that the panel is updated with
  \TRSM tasks, and the trailing sub-matrix is updated with \GEMM tasks. This
  avoids the row pivoting between the nodes usually performed by the
  classical LU factorization algorithm with partial pivoting, or by tournament
  pivoting algorithms~\cite{CALU}.
  
\paragraph*{b) \textsc{QR Step}} If the numerical criterion is not met, a QR
factorization has to be performed. Many solutions could be used for this
elimination step. We chose to exploit the \HQR method implementation
presented in~\cite{dongarra2013hierarchical}. This allowed us to experiment
with different kinds of reduction trees, so as to find the most adapted
solution to our problem. Our default tree (which we use in all of our
experiments) is a hierarchical tree made of \Greedy reduction trees inside
nodes and a \MC reduction tree between the nodes. The \MC tree between nodes has
been chosen for its short critical path and 
its good pipelining of consecutive trees, in case some QR steps are performed in
sequence. The \Greedy reduction tree is favored within a node. A two-level
hierarchical approach is natural when considering multicore parallel
distributed architectures. (See~\cite{dongarra2013hierarchical} for more
details on the reduction trees).

To implement the \luqralgo within the \parsec framework, two
extensions had to be implemented within the runtime. The first extension 
allows the programmer to generate data during the
execution with the \textsc{OUTPUT} keywords. This data is then inserted
into the tracking system of the runtime to follow its path in the
dataflow. This is what has been used to generate the backup on the fly, 
and to limit the memory peak of the algorithm.
A second extension has been made for the end detection of the
algorithm. Due to its distributed nature, \parsec goes over all the
domain space of each type of task of the algorithm 
and uses a predicate, namely \emph{the owner computes} rule, to decide if a task
is local or not. Local tasks are counted and the end of the
algorithm is detected when all of them have been executed.
Due to the dynamism in the \luqralgo, the size of the domain space is
larger than the number of tasks that will actually be executed. Thus, a
function to dynamically increase/decrease  the number of local tasks
has been added, so that the \emph{Propagate} task of each node updates the
local counter according to the elimination step chosen.

\section{Experiments}
\label{sec.experiments}

The purpose of this section is to present numerical experiments for the hybrid \luqralgo, and to highlight the trade-offs between stability and performance that can be achieved by tuning the threshold $\alpha$ in the robustness criterion (see 
Section~\ref{sec.criteria}). 

\subsection{Experimental framework}

We used \textit{Dancer}, a parallel machine hosted at the Innovative
Computing Laboratory (ICL) in Knoxville, to run the experiments.  This
cluster has 16 multi-core nodes, each equipped with
$8$ cores, and an Infiniband 10G interconnection network. The nodes
feature two Intel Westmere-EP E5606 CPUs at 2.13GHz. The system is running
the Linux 64bit operating system, version 3.7.2-x86\_64. The software
was compiled with the Intel Compiler Suite 2013.3.163. BLAS kernels were
provided by the MKL library and OpenMPI 1.4.3 has been used for the
MPI communications by the \parsec runtime. Each computational thread
is bound to a single core using the HwLoc 1.7.1 library. If not
mentioned otherwise, we will use all $16$ nodes and the
data will be distributed according to a $4$-by-$4$
2D-block-cyclic distribution. The theoretical peak performance of the 16 nodes is
1091 G\glops.

For each experiment, we consider a square tiled-matrix $A$ of size
$N$-by-$N$, where $N = n \times n_{b}$.  The tile size $n_b$ has been fixed
to $240$ for the whole experiment set, because this value was found to
achieve good performance for both LU and QR steps. We evaluate the
backward stability by computing the HPL3 accuracy test of the
High-Performance Linpack benchmark~\cite{cpe03}:
$$HPL3 = \frac{\|Ax - b\|_{\infty}}{\|A\|_{\infty}\|x\|_{\infty}
  \times \epsilon \times N},$$ 
  where $x$ is the computed solution and
$\epsilon$ is the machine precision. Each test is run with double
precision arithmetic.  For performance, we point out that the number
of floating point operations executed by the hybrid algorithm depends
on the number of LU and QR steps performed during the
factorization. Thus, for a fair comparison, we assess the efficiency
by reporting the \emph{normalized} G\glops performance computed as
$$\text{G\glops} = \frac{\frac{2}{3}N^3}{\textsc{Execution time}},$$
where $\frac{2}{3} N^3$ is the number of \flop for LU with partial
pivoting and $\textsc{Execution time}$ is the execution time of the
algorithm. With this formula, QR factorization will only achieve half of the
performance due to the $\frac{4}{3} N^3$ \flop of the algorithm.

\begin{figure*}[ht!]
\centering
\includegraphics[width=\linewidth]{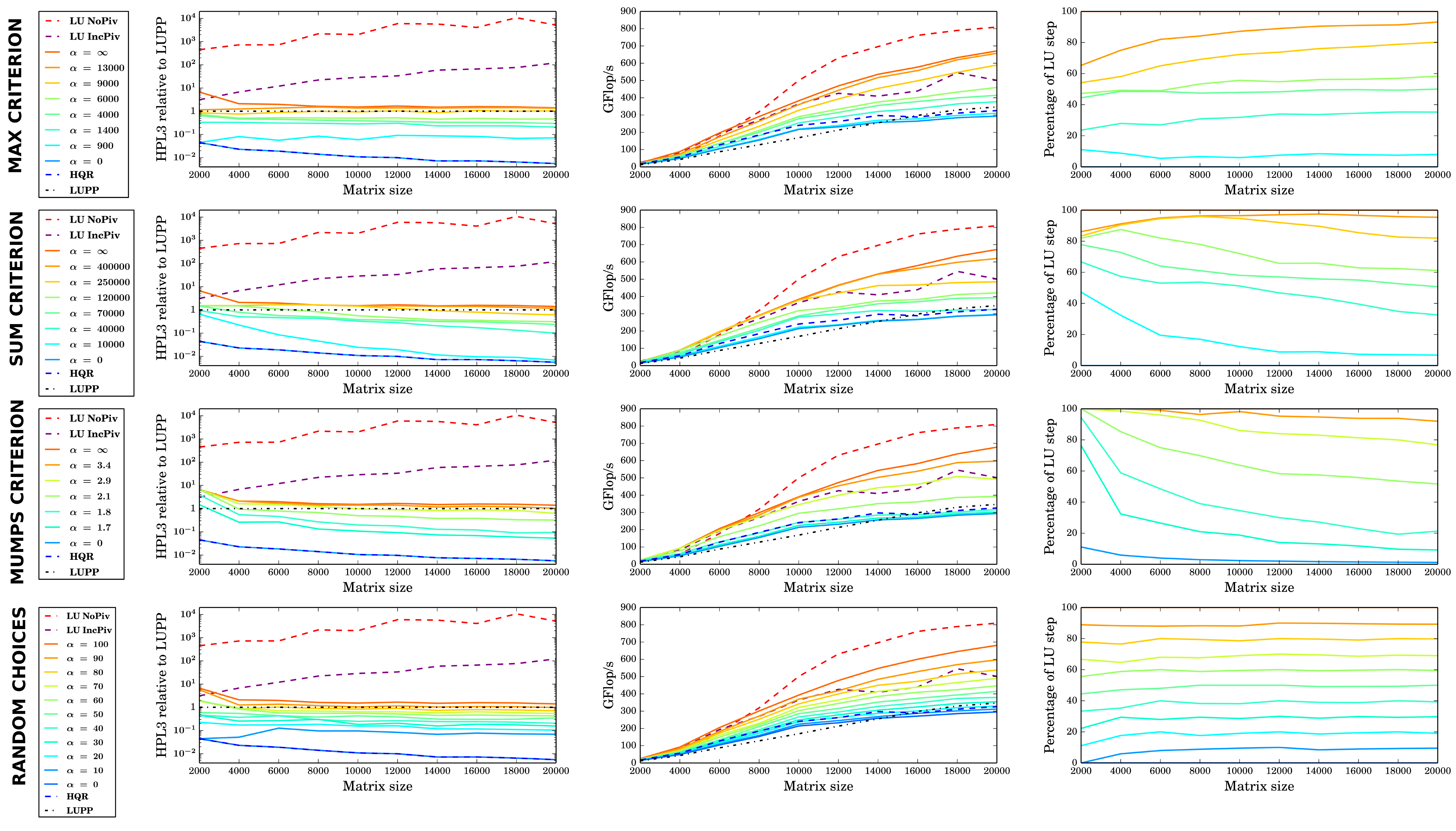}
\caption{Stability, performance, and percentage of LU steps obtained by the three criteria and by random choices, for random matrices, on the Dancer  platform (4x4 grid).}
\label{fig.random_mat}
\end{figure*}

\subsection{Results for random matrices}

We start with the list of the algorithms used for comparison with the \luqralgo. All of these methods are all implemented within the \dague framework:
\begin{compactitem}
\item LU NoPiv, which performs pivoting only inside the diagonal tile but no pivoting across tiles (known to be
both efficient and unstable)
\item LU IncPiv, which performs incremental pairwise pivoting across all tiles 
in the elimination panel~\cite{tileplasma,QuintanaOrti:2009} (still efficient but not stable either)
\item Several instances of the hybrid \luqralgo, for different values of the robustness parameter $\alpha$. Recall that the algorithm
performs pivoting only across the diagonal domain, hence involving no remote communication nor synchronization.
\item HQR, the Hierarchical QR
  factorization~\cite{dongarra2013hierarchical}, with the same
  configuration as in the QR steps of the \luqralgo:  \Greedy reduction trees inside nodes
  and \MC reduction trees between the nodes.
\end{compactitem}
For reference, we also include a comparison with PDGEQRF: this is the LUPP algorithm (LU
with partial pivoting across all tiles of the elimination panel) from the reference ScaLAPACK implementation~\cite{scalapack}.

\begin{table*}[htbp]
\begin{center}
\begin{tabular}{|c|c|c|c|c|c|c|c|c|}
\hline
	Algorithm & $\alpha$ & Time & \% LU steps & Fake G\glops & True G\glops & Fake \% Peak Perf. & True \% Peak Perf. \\
\hline
	LU NoPiv &				& 6.29	& 100.0	& 848.6	& 848.6	& 77.8	& 77.8 \\
	LU IncPiv&				& 9.25	& 100.0	& 576.4	& 576.4	& 52.9	& 52.9 \\
	LUQR (MAX) &	$\infty$	& 7.87	& 100.0	& 677.7	& 677.7	& 62.1	& 62.1 \\
	LUQR (MAX) &	13000		& 7.99	& 94.1	& 667.7	& 707.4	& 61.2	& 64.9 \\
	LUQR (MAX) &	9000		& 8.62	& 83.3	& 619.0	& 722.2	& 56.8	& 66.2 \\
	LUQR (MAX) &	6000		& 10.95	& 61.9	& 486.9	& 672.4	& 44.6	& 61.7 \\
	LUQR (MAX) &	4000		& 12.43	& 51.2	& 429.0	& 638.4	& 39.3	& 58.5 \\
	LUQR (MAX) &	1400		& 13.76	& 35.7	& 387.6	& 636.9	& 35.5	& 58.4 \\
	LUQR (MAX) &	900			& 16.39	& 11.9	& 325.4	& 612.0	& 29.8	& 56.1 \\
	LUQR (MAX) & 	0			& 18.05	& 0.0	& 295.5	& 590.9	& 27.1	& 54.2 \\ 
	HQR		 &				& 16.01	& 0.0	& 333.1	& 666.1	& 30.5	& 61.1 \\
	LUPP	 &				& 15.30	& 100.0	& 348.6	& 348.6	& 32.0	& 32.0 \\
\hline
\end{tabular}
\end{center}
\caption{Performance obtained by each algorithm, for $N=20,000$, on the Dancer platform 
	($4\times4$ grid). We only show the results for the \luqralgo with the Max criterion.  The other criteria have
	similar performance. In column Fake G\glops, we assume all algorithms perform $\frac{2}{3} N^3$ \flop.  
	In column True G\glops, we compute the number of \flop to be $( \frac{2}{3}f_{LU} + \frac{4}{3}(1-f_{LU}) ) N^3$, where $f_{LU}$ is the
	fraction of the steps that are LU steps (column 4).}\label{tab.perf} 
\end{table*}

Figure~\ref{fig.random_mat} summarizes all results for random matrices. 
It is organized as follows: each of the first three rows corresponds to one criterion. Within a row:
\begin{compactitem}
\item the first column shows the relative stability (ratio of HPL3 value divided by HPL3 value for LUPP)
\item the second column shows the \text{G\glops} performance
\item the third column shows the percentage of LU steps during execution
\end{compactitem}
The fourth row corresponds to a random choice between LU and QR at each step, and is intended to assess the performance obtained for a given ratio of LU vs QR steps. 
Plotted results are average values obtained on a set of 100 random matrices (we observe 
a very small standard deviation, less than $2\%$). 

For each criterion, we experimentally chose a set of values of $\alpha$ that
provides a representative range of ratios for the number of LU and QR steps. As
explained in Section~\ref{sec.criteria}, for each criterion, the smaller the
$\alpha$ is, the tighter the stability requirement. Thus, the numerical
criterion is met less frequently and the hybrid algorithm processes fewer LU
steps.  A current limitation of our approach is that we do not know how to
auto-tune the best range of values for $\alpha$, which seems to depend heavily
upon matrix size and available degree of parallelism. In addition, the
range of useful $\alpha$ values is quite different for each criterion.

For random matrices, we observe in Figure~\ref{fig.random_mat}
that the stability of LU~NoPiv and LU~IncPiv is not satisfactory. 
We also observe  that, for each criterion, small
values of $\alpha$ result in better stability, to the detriment of performance.
For $\alpha = 0$,  \luqralgo processes only QR steps, which leads to the exact
same stability as the HQR Algorithm and almost the same performance results.
The difference between the performance of \luqralgo with $\alpha = 0$ and HQR
comes from the cost of the decision making process steps (saving the panel,
computing the LU factorization with partial pivoting on the diagonal domain,
computing the choice, and restoring the panel). Figure~\ref{fig.random_mat}
shows that the overhead due to the decision making process is approximately
equal to 10\% for the three criteria. This overhead, computed when QR
eliminations are performed at each step, is primarily due to the
backup/restore steps added to the critical path when QR is
chosen. Performance impact of the criterion computation itself is
negligible, as one can see by comparing performance of the random
criterion to the MUMPS and Max criteria. 

\luqralgo with $\alpha = \infty$ and LU NoPiv both process only LU steps.  The
only difference between both algorithms in term of error analysis is that LU
NoPiv seeks for a pivot in the diagonal tile, while \luqralgo with $\alpha =
\infty$  seeks for a pivot in the diagonal domain.  This difference has a
considerable impact in term of stability, in particular on random matrices.
\luqralgo with $\alpha = \infty$  has a stability slightly inferior to that of
LUPP and significantly better to that of LU NoPiv.  When the matrix size
increases, the relative stability results of the \luqralgo with $\alpha =
\infty$ tends to $1$, which means that, on random matrices, processing an LU
factorization with partial pivoting on a diagonal domain followed by a direct
elimination without pivoting for the rest of the panel is almost as stable as
an LU factorization with partial pivoting on the whole panel.  A hand-waving
explanantion would go as follows.  The main instabilities are proportional to
the small pivots encountered during a factorization.  Using diagonal pivoting,
as the factorization of the diagonal tile proceeds, one is left with fewer and
fewer choices for a pivot in the tile.  Ultimately, for the last entry of the tile in
position ($n_b$,$n_b$), one is left with no choice at all.  When working on
random matrices, after having performed several successive diagonal
factorizations, one is bound to have encountered a few small pivots.  These
small pivots lead to a bad stability.  Using a domain (made of several tiles) for the factorization
significantly increases the number of choice for the pivot and it is not any
longer likely to encounter a small pivot. Consequently diagonal domain pivoting
significantly increases the stability of the \luqralgo with $\alpha = \infty$.
When the local domain gets large enough (while being significanty less than
$N$), the stability obtained on random matrices is about the same as partial
pivoting. 

When $\alpha = \infty$, our criterion is deactivated and our algorithm always performs LU step.
We note that, when $\alpha$ is reasonable, (as opposed to $\alpha = \infty$,)
the algorithm is stable whether we use a diagonal domain or a diagonal tile.
However using a diagonal domain increases the chance of well-behaved pivot tile
for the elimination, therefore using a diagonal domain (as opposed to a
diagonal tile) increases the chances of an LU step.

Using random choices leads to results comparable to those obtained with the three criteria. 
However, since we are using random matrices in this experiment set, we need to be careful before drawing any conclusion
on the stability of our algorithms. If an algorithm is not stable on random matrices, 
this is clearly bad. 
However we cannot draw any definitive conclusion if an algorithm is stable for random matrices.

Table~\ref{tab.perf} displays detailed performance results for the Max criterion with $N=20,000$.
In column Fake G\glops, we assume all algorithms perform $\frac{2}{3} N^3$ \flop.  
In column True G\glops, we compute the number of \flop to be 
$( \frac{2}{3}f_{LU} + \frac{4}{3}(1-f_{LU}) ) N^3$, where $f_{LU}$ is the
fraction of the steps that are LU steps (column 4).  For this example, we see that the \luqralgo
reaches a peak performance of 677.7 G\glops (62.1\% of the theoretical peak) when every
step is an LU step.  Comparing \luqralgo with $\alpha=0$ and 
HQR shows the overhead for the decision making process is 12.7\%.  
We would like the true performance to be constant as the number of QR steps increases.  For this
example, we see only a slight decrease in performance, from 62.1\% for $\alpha=\infty$ to 54.2\% for $\alpha=0$.
HQR maintains nearly the same true performance (61.1\%) as the \luqralgo with $\alpha=\infty$.
Therefore, the decrease in true performance is due largely to the overhead of restoring the panel.

\subsection{Results for special matrices}
\label{sec.stability}

For random matrices, we obtain a good stability with random choices,  almost as good as with the three criteria.
However, as mentioned above, we should draw no definite conclusion. 
To highlight the need for a smart criterion, we tested the hybrid \luqralgo on a collection of
matrices that includes several pathological matrices on which LUPP fails because of large growth factors. This
set of special matrices described in Table~\ref{tab.spec} includes
ill-conditioned matrices as well as sparse matrices, and mostly comes from the
Higham's Matrix Computation Toolbox~\cite{Higham2002}.

\begin{table*}
\small{
\begin{tabular}{rll}
\textbf{No.} & \textbf{Matrix} & \textbf{Description} \\
1 & house & Householder matrix, $A = eye(n) - \beta * v * v$ \\
2 & parter & Parter matrix, a Toeplitz matrix with most of singular values near $\Pi$. $A(i, j) = 1/(i - j + 0.5)$. \\
3 & ris & Ris matrix, matrix with elements $A(i, j) = 0.5/(n - i - j + 1.5)$. The
eigenvalues cluster around $-\Pi/2$ and $\Pi/2$. \\
4 & condex & Counter-example matrix to condition estimators. \\
5 & circul & Circulant matrix \\
6 & hankel & Hankel matrix, $A = hankel(c, r)$, where $c=randn(n, 1)$, $r=randn(n, 1)$, and $c(n) = r(1)$. \\
7 & compan & Companion matrix (sparse), $A = compan(randn(n+1,1))$. \\
8 & lehmer & Lehmer matrix, a symmetric positive definite matrix such that
$A(i, j) = i/j$ for $j \geq i$. Its inverse is tridiagonal. \\
9 & dorr & Dorr matrix, a diagonally dominant, ill-conditioned, tridiagonal matrix
(sparse). \\
10 & demmel & $A = D*(eye(n) + 10-7 *rand(n))$, where $D = diag(1014*(0:n-1)/n )$. \\
11 & chebvand & Chebyshev Vandermonde matrix based on n equally spaced points on
the interval [0,1]. \\
12 & invhess & Its inverse is an upper Hessenberg matrix. \\
13 & prolate & Prolate matrix, an ill-conditioned Toeplitz matrix. \\
14 & cauchy & Cauchy matrix. \\
15 & hilb & Hilbert matrix with elements $1/(i + j - 1)$. $A =hilb(n)$. \\
16 & lotkin & Lotkin matrix, the Hilbert matrix with its first row altered to all ones. \\
17 & kahan & Kahan matrix, an upper trapezoidal matrix. \\
18 & orthogo & Symmetric eigenvector matrix: $A(i,j) = sqrt(2/(n+1)) * sin(i*j*\pi/(n+1))$\\
19 & wilkinson & Matrix attaining the upper bound of the growth factor of GEPP. \\
20 & foster & Matrix arising from using the quadrature method to solve a certain Volterra integral equation. \\
21 & wright & Matrix with an exponential growth factor when Gaussian elimination with Partial Pivoting is used. \\
\end{tabular}	
}
\caption{\label{tab.spec} Special matrices in the experiment set.}
\end{table*}

\begin{figure}[htbp]
\centering
\includegraphics[width=\linewidth]{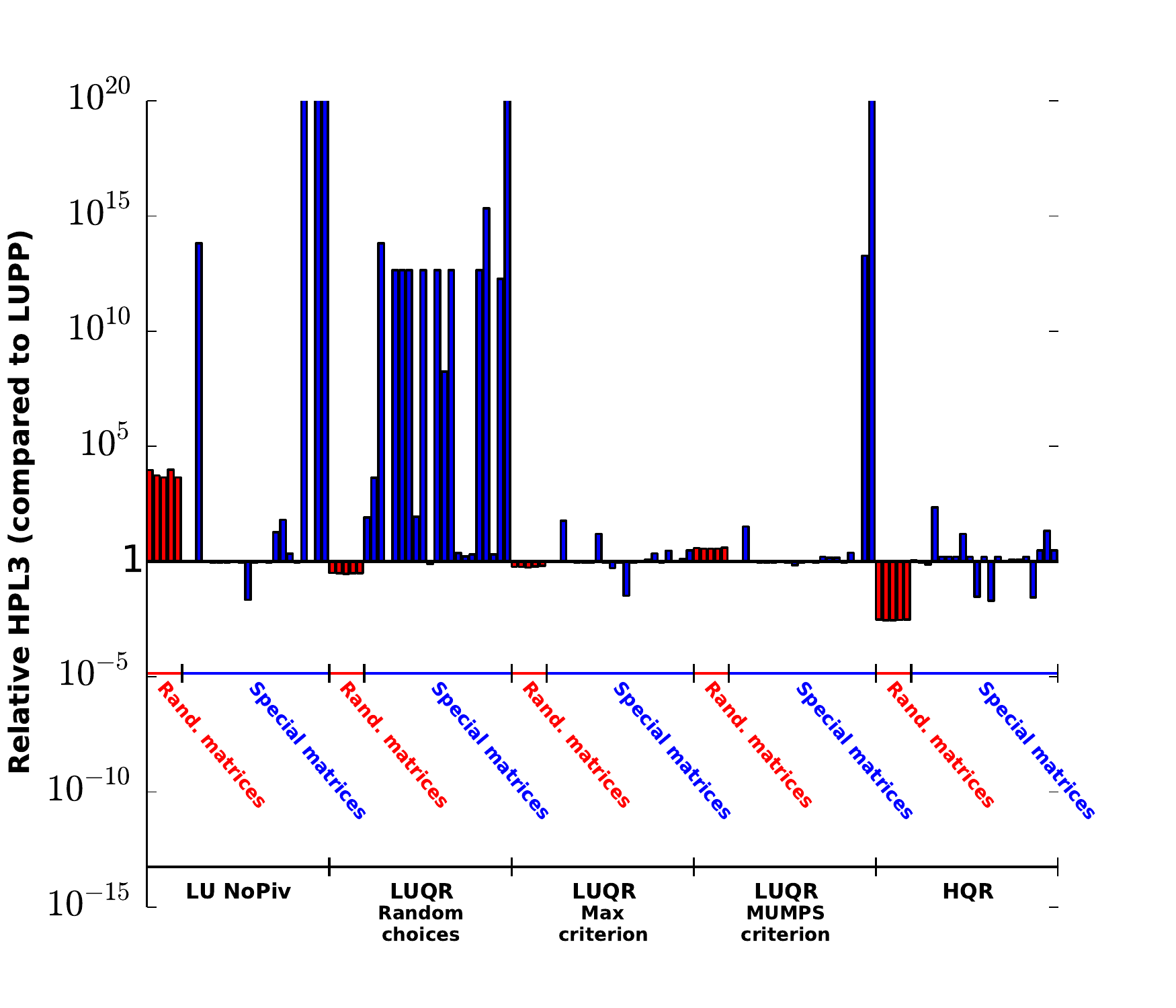}
\caption{Stability on special matrices.}
\label{fig.spec}
\end{figure}

Figure~\ref{fig.spec} provides the relative stability (ratio of HPL3 divided by HPL3 for LUPP) 
obtained by running the hybrid \luqralgo on a set of $5$ random
matrices and on the set of special matrices. Matrix size is set to
$N=40,000$, and experiments were run on a $16$-by-$1$ process grid. The parameter $\alpha$ has been set to $50$ for the random criterion, $6,000$ for the Max criterion, and $2.1$ for the MUMPS criterion (we do not report result for the Sum criterion because they are the same as they are for Max).
Figure~\ref{fig.spec}  considers LU NoPiv, HQR and the \luqralgo . The first observation is that using random choices 
now leads to numerical instability. The Max criterion provides a good stability ratio on every tested matrix (up to 58 for the RIS matrix and down to 0.03 for the Invhess matrix). The MUMPS criterion also gives modest growth factor for the whole experiment set except for the Wilkinson and the Foster matrices, for which it fails to detect some ``bad'' steps. 

We point out that we also experimented with the Fiedler matrix from  Higham's Matrix Computation Toolbox~\cite{Higham2002}. We observed that LU NoPiv and LUPP failed (due to small values rounded up to 0 and then illegally used in a division), while the Max and the MUMPS criteria provide HPL3 values ($\approx 5.16\times10^{-09}$ and $\approx 2.59\times10^{-09}$) comparable to that of HQR ($\approx 5.56\times10^{-09}$). This proves that our criteria can detect and handle pathological cases for which the generic LUPP algorithm fails.

\subsection{Assessment of the three criteria}

With respect to stability, while the three criteria behave similarly on random matrices, we observe different behaviors for special matrices.
The MUMPS criterion provides good results for most of the tested matrices but not for all. If stability is the key concern,
one may prefer to use the Max criterion (or the Sum criterion), which performs well for all special matrices (which means that the upper bound of $(1 + \alpha) ^{n-1}$ on the growth is quite pessimistic).

With respect to performance, we observe very comparable results, which means that the overhead induced
by computing the criterion at each step is of the same order of magnitude for all criteria.

The overall conclusion is that all criteria bring significant improvement over LUPP in terms of stability, and
over HQR in terms of performance. Tuning the value of the robustness parameter $\alpha$ enables the exploration of
a wide range of stability/performance trade-offs.

\section{Related work}
\label{sec.related}

State-of-the-art QR factorizations use multiple eliminators per panel, 
in order to dramatically reduce the critical path of the algorithm. These algorithms
are unconditionally stable, and their parallelization has been fairly well studied
on shared memory systems~\cite{Buttari2008,QuintanaOrti:2009,sc-paper2011} and
on parallel distributed systems~\cite{dongarra2013hierarchical}. 

The idea of mixing Gaussian transformations and orthogonal transformations has been considered once before.
Irony and Toledo~\cite{irony2006snap} present an algorithm for reducing a banded symmetric indefinite matrix
to diagonal form.  The algorithm uses symmetric Gaussian transformations and Givens rotations to maintain the
banded symmetric structure and maintain similar stability to partial symmetric pivoting.  

The reason for using LU kernels instead of
QR kernels is performance: (i) LU performs half the number of \flop of QR; (ii) LU
kernels relies on GEMM kernels which are very efficient while QR kernels are more complex and much
less tuned, hence not that efficient; and (iii) the LU update is much more parallel
than the QR update. So all in all, LU is much faster 
than QR (as observed in the performance results of Section~\ref{sec.experiments}).
Because of the large number of communications and synchronizations induced by pivoting in the reference LUPP algorithm,
\emph{communication-avoiding} variants of LUPP have been introduced~\cite{dghl:sisc:12}, 
but they have proven much more challenging to design because of stability issues.
In the following, we review several approaches:

\subsection{LUPP}
LU  with partial pivoting is not a communication-avoiding scheme and
its performance in a parallel distributed environment is low (see Section~\ref{sec.experiments}).
However, the LUPP algorithm is {\em stable in practice}, and we use it as a reference for stability.

\subsection{LU~NoPiv}
The most basic communication-avoiding LU algorithm is LU~NoPiv. This
algorithm is stable for block diagonal dominant
matrices~\cite{Higham2002,blocklu:95}, but breaks down if it encounters a nearly
singular diagonal tile, or loses stability if it encounters  a diagonal tile whose
smallest singular value is too small.

Baboulin et al.~\cite{Baboulin:2012:ALS} propose to apply  a random
transformation to the initial matrix, in order to use LU~NoPiv while maintaining
stability.  This approach gives about the same performance as LU~NoPiv, since
preprocessing and postprocessing costs are negligible. It is hard to be
satisfied with this approach~\cite{Baboulin:2012:ALS} because for any matrix
which is rendered stable by this approach (i.e, LU~NoPiv is stable), there
exists a matrix which is rendered not stable. Nevertheless, in practice, this
proves to be a valid approach.

\subsection{LU~IncPiv}
LU~IncPiv is another communication-avoiding LU
algorithm~\cite{tileplasma,QuintanaOrti:2009}. Incremental pivoting is also
called {\em pairwise pivoting}. The stability of the algorithm~\cite{tileplasma} is
not sufficient and degrades as the number of tiles in the matrix increases
(see our experimental results on random matrices). The method also suffers some
of the same performance degradation of QR factorizations with multiple eliminators per panel,  namely
low-performing kernels, and some dependencies in the update phase.

\subsection{CALU}
CALU~\cite{CALU} is a communication-avoiding LU. It uses tournament pivoting which has been proven
to be {\em stable in practice}~\cite{CALU}.  
CALU shares the (good) properties of one of our LU steps: 
(i) low number of \flop; (ii) use of efficient GEMM kernels; and (iii) embarrassingly parallel update.
The advantage of CALU over our algorithm is essentially that it performs only LU steps, while
our algorithm might need to perform some (more expensive) QR steps.
The disadvantage is that, at each step, CALU needs to perform global pivoting on the whole panel,
 which then needs to be reported
during the update phase to the whole trailing submatrix.
There is no publicly available implementation of 
parallel distributed CALU, and it was not possible to compare stability or performance.
CALU is known to be stable in practice~\cite{grigori2008calu,lawn280}.
Performance results of CALU in parallel distributed are presented in~\cite{grigori2008calu}.
Performance results of CALU on a single multicore node are presented in~\cite{lawn280}.

\section{Conclusion}
\label{sec.conclusion}

Linear algebra software designers have been struggling for years to improve the
parallel efficiency of LUPP (LU with partial pivoting), the de-facto choice
method for solving dense systems. The search for good pivots throughout the
elimination panel is the key for stability (and indeed both NoPiv and IncPiv
fail to provide acceptable stability), but it induces 
several short-length communications that dramatically decrease the overall
performance of the factorization.

Communication-avoiding algorithms are a recent alternative which proves very
relevant on today's archictectures.  For example, in our experiments, our HQR
factorization~\cite{dongarra2013hierarchical} based of QR kernels ends
with similar performance as ScaLAPACK LUPP while performing 2x more \flop, using slower
sequential kernels, and a less parallel update phase. In this paper,
stemming from the key observation that LU steps and QR steps can be mixed
during a factorization, we present the \luqralgo whose goal is to accelerate
the HQR algorithm by introducing some LU steps whenever these do not compromise
stability.
The hybrid algorithm represents
dramatic progress in a long-standing research problem. By restricting to
pivoting inside the diagonal domain, i.e., locally, but by doing so only when
the robustness criterion forecasts that it is safe (and going to a QR step
otherwise), we improve performance while guaranteeing stability. And we provide
a continuous range of trade-offs between LU NoPiv (efficient but only stable
for diagonally-dominant matrices) and QR (always stable but twice as costly and
with less performance).
For some classes of matrices (e.g., tile diagonally dominant),
the \luqralgo will only perform LU steps.

 This work opens several research directions. First, as already mentioned, the
 choice of the robustness parameter $\alpha$ is left to the user, and it would
 be very interesting to be able to auto-tune a possible range of values as a
 function of the problem and platform parameters. Second, there are many
 variants and extensions of the hybrid algorithm that can be envisioned. Several
 have been mentioned in Section~\ref{sec.luqr-algo}, and many others could be
 tried. Another goal would be to derive LU 
 algorithms with several eliminators per panel (just as for HQR) to decrease the critical path,
provided the availability of a reliable robustness test to ensure stability.

\subsection*{Acknowledgment}

The research of Julien Langou and Bradley R. Lowery was fully supported by the National Science Foundation grant \# NSF CCF 1054864. The research of Yves Robert was supported in part by the French ANR \emph{Rescue} project and by the 
Department of Energy \# DOE DE-SC0010682. Yves Robert is with Institut Universitaire de France. The authors thank Thomas
Herault for his help with \dague, and Jean-Yves L'Excellent for discussions on stability estimators within the MUMPS software.

\bibliographystyle{IEEEtran}
\bibliography{biblio}

\end{document}